\documentclass[11pt]{amsart}

\usepackage{graphicx}
\usepackage{amssymb,mathrsfs}
\usepackage{subfigure}

\usepackage{enumerate}        

\newcommand {\bng}{B_n\Gamma}
\newcommand {\pbng}{PB_n\Gamma}
\newcommand {\btg}{B_2\Gamma}
\newcommand {\cng}{\mathcal{C}^n\Gamma}
\newcommand {\ucng}{U\mathcal{C}^n\Gamma}

\newcommand {\dng}{\mathcal{D}^n\Gamma}
\newcommand {\udng}{U\mathcal{D}^n\Gamma}

\newcommand {\ud}[2]{U\mathcal{D}^{#1} #2}

\newtheorem{theorem}{Theorem}[section]

\newtheorem{proposition}[theorem]{Proposition}
\newtheorem{corollary}[theorem]{Corollary}

\theoremstyle{definition}
\newtheorem{definition}[theorem]{Definition}
\newtheorem{convention}[theorem]{Convention}
\newtheorem{example}[theorem]{Example}

\begin{document}

\title[Graph groups and graph braid groups]{Embedding Right-Angled Artin
Groups into
Graph Braid Groups}
\author[L.Sabalka]{Lucas~Sabalka$\!\,^1$}
      \address{\tt Department of Mathematics\\
               U. of Illinois at Urbana-Champaign\\
               Champaign IL 61820
\newline       http://www.math.uiuc.edu/\~{}sabalka}
      \email{sabalka at math.uiuc.edu}

\footnotetext[1]{Research supported by an NSF Graduate Research
Fellowship and a VIGRE graduate fellowship from the UIUC
mathematics department.}

\begin{abstract}
We construct an embedding of any right-angled Artin group $G(\Delta)$ 
defined by a graph $\Delta$ into a graph braid group.  The number of 
strands required for the braid group is equal to the chromatic number of 
$\Delta$.  This construction yields an example of a hyperbolic surface 
subgroup embedded in a two strand planar graph braid group.
\end{abstract}

\maketitle

\section{Introduction}\label{sec:intro}

Let $\Delta$ be a finite simple graph.  We examine a group $A =
G(\Delta)$, called the \emph{right-angled Artin group} associated to
$\Delta$, defined with the following presentation:  for each vertex $a_i$
of $\Delta$, there exists a corresponding generator for $A$, and two
generators $a_i$ and $a_j$ commute (for $i \neq j$) if and only if the
vertices $a_i$ and $a_j$ are connected by an edge in $\Delta$.  In
particular, a right-angled Artin group is a group which has a presentation
where the only relators are commutators of generators.  Right-angled Artin
groups are of interest as sources of subgroups with complicated homology
and homotopy finiteness properties; see for example \cite{Stallings}.  
Right-angled Artin groups are well studied; for a more extensive
reference, see for example \cite{BestvinaBrady},
\cite{DavisJanuszkiewicz}.

Given a graph $\Gamma$, the \emph{ordered configuration space} of $n$
points on $\Gamma$, denoted $\cng$, is the open subset of the direct
product $\Pi^n\Gamma$ obtained by removing the \emph{diagonal}, $Diag :=
\{(x_1,\dots,x_n) \in \Pi^n\Gamma | x_i = x_j \hbox{ for some } i \neq
j\}$.  The fundamental group of the ordered configuration space $\cng$ of
$n$ points on some graph $\Gamma$, denoted $\pbng$, is called a \emph{pure
graph braid group}.  The \emph{unordered configuration space} of $n$
points on $\Gamma$, denoted $\ucng$, is the quotient of $\cng$ by the
action of the symmetric group permuting the factors.  Thus, $\ucng :=
\left\{\{x_1,\dots,x_n\} | x_i \in \Gamma, x_i = x_j \hbox{ for some } i
\neq j\right\}$, where we use braces to indicate that the coordinates are
unordered.  The fundamental group $\pi_1\ucng$ of $n$ strands on $\Gamma$
is called a \emph{graph braid group}, denoted $\bng$.

Graph braid groups are of interest because of their connections with
robotics and mechanical engineering.  Graph braid groups can, for
instance, model the motions of robots moving about a factory floor
(\cite{Ghrist}, \cite{Farber1}, \cite{Farber2}), or the motions of
microscopic balls of liquid on a nano-scale electronic circuit
(\cite{GhristPeterson}). For more information about graph braid groups see
for instance \cite{Abrams1}, \cite{FarleySabalka1}, \cite{AbramsGhrist}.

In \cite{Ghrist}, Ghrist conjectured that every graph braid group is a
right-angled Artin group.  In \cite{Abrams1}, Abrams was able to find
counterexamples to this conjecture: namely, $PB_2K_5$ and $PB_2K_{3,3}$
are surface groups.  Connolly and Doig \cite{ConnollyDoig} proved a
partial positive result, that braid groups on \emph{linear trees} are
right-angled Artin, where a tree is linear if there exists an embedded
interval containing all of the vertices of degree at least 3. In general,
though, it appears that graph braid groups are usually not right-angled
Artin groups.  Mautner reported via private communication that he has some
examples of graph braid groups on planar graphs with nontrivial
fundamental group which are not right-angled Artin \cite{Mautner}.  More
recently, in \cite{FarleySabalka2a} it is proved that tree braid groups 
are right-angled Artin if and only if the tree is linear or there are less
than 4 strands.

A related line of inquiry asks whether there exist embeddings between
right-angled Artin groups and graph braid groups.  A positive result in
this direction is due to Crisp and Wiest:  for any finite graph $\Gamma$ 
and any $n$, there exists a graph $\Delta$ such that $\bng$ embeds in 
$G(\Delta)$.

The goal of this paper is to prove the opposite direction of Crisp and 
Wiest's theorem:
\begin{theorem}\label{thm:embedding}
For every finite graph $\Delta$ and any coloring $C$ of $\Delta$ with $n$ 
colors, there exists a graph $\Gamma$ such that the right-angled Artin 
group $G(\Delta)$ embeds into the graph braid group $\bng$.
\end{theorem}

For example, letting $C$ be the trivial coloring assigning to each vertex
of $\Delta$ a different color, we obtain an embedding of $G(\Delta)$ into
$\bng$ for some graph $\Gamma$, where $n$ is the number of vertices of
$\Delta$.  At the other extreme, $n$ can be as small as the chromatic
number of $\Delta$, by letting $C$ be a coloring realizing the chromatic
number.

In Section \ref{sec:embeddings}, we use Theorem \ref{thm:embedding} to
obtain an explicit embedding of the right-angled Artin group associated to
the cycle of length 6, $C_6$, into a two strand planar graph braid group.  
By a result of Servatius, Droms, and Servatius
\cite{ServatiusDromsServatius}, the group $G(C_6)$ itself contains an
embedded hyperbolic surface subgroup (where a \emph{hyperbolic surface
group} is the fundamental group of a compact hyperbolic surface without
boundary).  Thus, we prove as a corollary the following theorem:

\begin{theorem}\label{thm:planartwostrandsurface}
There exists a \emph{planar} graph braid group which contains a hyperbolic 
surface subgroup, with only $n = 2$ strands.
\end{theorem}

There were previously no known examples of hypberbolic surface subgroups 
in planar graph braid groups.

The proof presented here of Proposition \ref{prop:psiinjective} was
discovered by Dan Farley and is a refinement of the original proof of the
author's.  The current construction of the graph $\Gamma$ in Section
\ref{sec:preliminaries} was inadvertently suggested by Bert Weist and is a 
slight improvement over the original construction.  REU student Go Fujita
recognized (also inadvertently) that the number of strands used in 
Proposition \ref{prop:psiinjective} corresponds to a coloring of $\Delta$.

The author would like to thank Dan Farley for numerous helpful discussions
on this and related material, as well as his advisor, Ilya Kapovich.  The
author is also thankful for the (accidental) helpful observations of
Bert Wiest and Go Fujita made while trying to understand an earlier proof.

\section{Preliminaries}\label{sec:preliminaries}

A \emph{coloring} $C$ of a finite simple graph $\Delta$ with $n$ colors is
a function from the vertices $V\Delta$ of $\Delta$ to the finite set $S =
\{1, \dots, n\}$ such that if two vertices $v$ and $w$ in $\Delta$ are
neighbors then $C(v) \neq C(w)$.  The elements of $S$ are called
\emph{colors}.  The \emph{color} of a vertex is its image under $C$.  The
\emph{chromatic number} of $\Delta$ is the smallest number of colors
needed to have a coloring of $\Delta$.

An \emph{edge loop} in a CW-complex is a closed path consisting of 
a sequence of $1$-cells in the space.  

\begin{definition}[Halo] \label{def:halo}
Let $\Delta$ be a finite simple graph and let $C$ be a coloring of
$\Delta$ with $n$ colors $\{1,\dots,n\}$.  A connected graph $\Gamma$ is
called a \emph{halo of $\Delta$}, or a \emph{$\Delta$-halo}, if for each
vertex $a_i$ in $\Delta$ (equivalently, for each generator in $A$) there 
exists a simple edge loop $\gamma_i$ in $\Gamma$, called an \emph{Artin 
loop} of $\Gamma$, such that:
\begin{itemize}
\item for each color $c \in \{1, \dots, n\}$, there exists a vertex $x_c 
\in \Gamma$ common to all Artin loops of vertices colored $c$ and to no 
other Artin loops.

\item if $a_i$ and $a_j$ are not connected by an edge in $\Delta$, then
the Artin loops $\gamma_i$ and $\gamma_j$ intersect in exactly one vertex
$v$.  If $C(a_i) = C(a_j)$ then $v = x_{C(a_i)}$; otherwise, $v$ is on no
other Artin loop.

\item if $a_i$ and $a_j$ are connected by an edge in $\Delta$ (i.e. if
$a_i$ and $a_j$ commute), then the Artin loops $\gamma_i$ and $\gamma_j$
are disjoint.
\end{itemize}
We call the point $\{x_1, \dots, x_n\} \in \ucng$ the \emph{Artin 
basepoint} of $\Gamma$.
\end{definition}

\begin{figure}[htp]
\begin{center}
\subfigure[$\Delta$]{
  \includegraphics[width=2in]{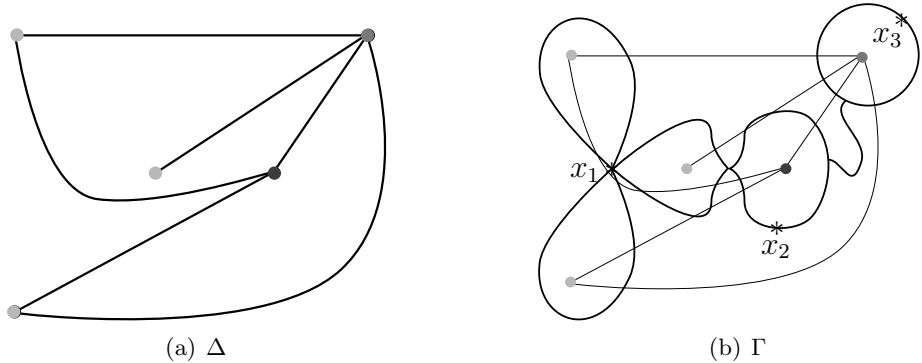}
  \label{fig:loopconstructionleft}
}
\hfill
\subfigure[$\Gamma$]{
  \input{loopconstructionb2.tex}
  \label{fig:loopconstructionright}
}
\end{center}

\caption{Constructing $\Delta$-halos.}{In \ref{fig:loopconstructionleft},
we have shown a possible graph $\Delta$ corresponding to some right-angled
Artin group along with a coloring of $\Delta$.  Then, in
\ref{fig:loopconstructionright}, we have shown a possible $\Delta$-halo
$\Gamma(\Delta)$, with the graph $\Delta$ subimposed in gray.}

\label{fig:loopconstruction}
\end{figure}

Although there is no clear canonical choice for a halo of an arbitrary
$\Delta$, when a particular $\Delta$-halo $\Gamma$ is specified we write
$\Gamma = \Gamma(\Delta)$. That $\Delta$-halos exist is 
clear; see for example Figure \ref{fig:loopconstruction}.

\begin{convention}\label{conv:Artinloops}
For any Artin basepoint of $\Gamma = \Gamma(\Delta)$, arbitrarily fix a
direction in which to traverse each $\gamma_i$.  Let $c$ be the color of
the vertex $a_i$.  We call the loop in $\ud{n}{\Gamma}$ from
$\{x_1,\dots,x_n\}$ corresponding to the strand at $x_c$ making a single
traversal of the Artin loop $\gamma_i$ in the chosen direction also by the
name of the Artin loop $\gamma_i$.
\end{convention}

The notion of a $\Delta$-halo is similar in flavor to the proof of Theorem 
10 of \cite{CrispWiest}, though the notions were discovered independently.  

Recall that the ordered configuration space $\cng$ is defined to be a
subset of the direct product $\Pi^n\Gamma$ of $n$ copies of $\Gamma$,
where we subtract off the diagonal $Diag := \{(x_1,\dots,x_n)|x_i = x_j
\hbox{ for some } i \neq j\}$.  The unordered configuration space $\ucng$
is then the quotient of this space by the action of $S_n$ permuting the
factors.  As $\Gamma$ is a graph, it is also a $1$-dimensional CW-complex,
where the interiors of edges are the open $1$-cells.  This means that
there is a product CW-complex structure on $\Pi^n\Gamma$, where an open
$k$-cell consists of an ordered collection of exactly $k$ edges and $n-k$
vertices of $\Gamma$.  Let $Diag'$ be the set of all open cells whose 
closure intersects $Diag$.  Thus
  $$Diag' = \{(y_1,\dots,y_n) | y_i \hbox{ is a cell in } \Gamma, 
  \hbox{ and } \overline{y_i} \cap \overline{y_j} \neq \emptyset 
  \hbox{ for some } i \neq j\}.$$
Let $\dng$ be the CW-complex $\Pi^n\Gamma - Diag'$, the 
\emph{discretized configuration space} of $n$ strands on $\Gamma$, and 
let $\udng$ be the quotient of $\dng$ given by permuting the factors of 
$\Pi^n\Gamma$ via the action of $S_n$, the \emph{unlabelled discretized 
configuration space} of $n$ strands on $\Gamma$.

For a general graph, call the vertices of valence at least 3
\emph{essential}.  In \cite{Abrams1}, Abrams proved that if $\Gamma$ is
sufficiently subdivided, then $\ucng$ is homotopy equivalent to $\udng$,
where sufficiently subdivided means every nontrivial path connecting two
essential vertices has at least $n - 1$ edges, and every nontrivial edge
loop has at least $n+1$ edges. In particular, this implies that $\bng
\cong \pi_1 \udng$ for $\Gamma$ sufficiently subdivided.  It is clear that
by subdividing edges of a graph $\Gamma$ we do not change the associated
graph braid group, so we always assume that our graph $\Gamma$ is
sufficiently subdivided.  Note the loops Artin $\gamma_i$ in $\ucng$ are
in fact simple edge loops in $\udng$ since the Artin basepoint $x = 
\{x_1,\dots,x_n\}$ is a $0$-cell in $\udng$.

\section{Embeddings and Corollaries}\label{sec:embeddings}

In \cite{Droms}, Droms proved that right-angled Artin groups have a
rigidity property:  two graphs $\Delta$ and $\Delta'$ are graph isomorphic
if and only if the groups $G(\Delta)$ and $G(\Delta')$ are group
isomorphic.  Thus, we may refer to a right-angled Artin group and the
graph which defines it interchangeably.

Let $\Delta$ be a simple graph with corresponding right-angled Artin group
$A = G(\Delta)$.  For any vertex $v \in \Delta$, let $link(v)$ denote the
full subgraph of $\Delta$ whose vertices are exactly those vertices
adjacent to $v$ in $\Delta$, not including $v$ itself.  Let $\Delta - v$
denote the subgraph of $\Delta$ formed by deleting the vertex $v$ and any
edge adjacent to $v$.  The group $A$ has the structure of an HNN-extension
\cite{BridsonHaefliger} with stable letter $v$, where the base group is
$G(\Delta - v)$, the associated subgroup is $G(link(v))$, and $v$
conjugates a given element of $G(link(v))$ to itself.

Britton's Lemma \cite{LyndonSchupp} tells us that if $g_0v^{\epsilon_1}g_1
v^{\epsilon_2}g_2\dots v^{\epsilon_k}g_k = 1$, where $\epsilon_i \in \{\pm
1\}$ and $g_i \in G(\Delta - v)$, then there must exist some $i$ such that
$\epsilon_{i+1} = - \epsilon_i$ and $g_i \in G(link(v))$.  Thus,
$v^{\epsilon_i} g_i v^{\epsilon_{i+1}} = g_i$, and it must be that $k$ is
not minimal.  We call $v^{\epsilon_i} g_i v^{\epsilon_{i+1}}$ a
\emph{$v$-pinch}.

We now recall the embedding of Crisp and Wiest of graph braid groups
into right-angled Artin groups, given by Theorem 2 of \cite{CrispWiest}.  
Let $\Gamma$ be a sufficiently subdivided graph.  Let $\Delta_{\Gamma}$ be
the graph whose vertex set is the edge set of $\Gamma$, and where two
vertices $e$ and $e'$ are connected by an edge if and only if the closures
of $e$ and $e'$ are disjoint in $\Gamma$.  Let $A_{\Gamma} :=
G(\Delta_{\Gamma})$ be the associated right-angled Artin group.

Crisp and Wiest prove their Theorem 2 by giving an embedding of
$\ud{n}{\Gamma}$ into an Eilenberg-MacLane complex for $A_{\Gamma}$ which
is a locally CAT(0) cubed complex.  We apply their argument on the level
of groups. In particular, let $\gamma$ be an edge loop of length $k$ in
$\ud{n}{\Gamma}$.  Each edge of $\gamma$ exactly corresponds to one strand
crossing one edge of $\Gamma$ while the other strands do not move.  Let
$e_1, \dots, e_k$ be the edges of $\Gamma$ corresponding to each edge of
$\gamma$ given in order. Let $w$ be the word in the generators of
$A_{\Gamma}$ given by $e_1\dots e_k$.  Let $\Phi: \bng =
\pi_1(\Gamma,\{x_1,\dots,x_n\}) \to A_{\Gamma}$ be defined by mapping the
homotopy class of $\gamma$ in $\pi_1 \ud{n}{\Gamma}$ to the element
represented by $w$ in $A_{\Gamma}$. Thus $\Phi$ is the ``forgetful" map,
ignoring all of the strands in $\gamma$ except the strand crossing an edge
in $\Gamma$. The map $\Phi$ is the Crisp and Wiest embedding.

The following proposition implies Theorem \ref{thm:embedding}.

\begin{proposition} \label{prop:psiinjective}
Let $\Delta$ be a finite graph and let $C$ be a coloring of $\Delta$ using
$n$ colors $\{1,\dots,n\}$.  Let $\Gamma = \Gamma(\Delta)$ be a halo of
$\Delta$ with Artin basepoint $\{x_1, \dots, x_n\}$, and let $a_i$ denote
the vertices of $\Delta$ for $i = 1, \dots$.  Set $A := G(\Delta)$.  The
map $a_i \mapsto \gamma_i^2$ (following Convention \ref{conv:Artinloops})
induces an injective homomorphism $\Psi: G(\Delta) \to \bng$.
\end{proposition}

We will see below that the square in the definition of $\Psi$ is necessary 
to obtain injectivity.  

\begin{proof} 
If $a_i$ and $a_j$ commute then so do $\Psi(a_i)$ and $\Psi(a_j)$, as
$\gamma_i$ and $\gamma_j$ are disjoint (Artin) loops in $\Gamma$ and the
associated loops in $\ud{n}{\Gamma}$ involve distinct strands.  Thus the
map $a_i \mapsto \gamma_i^2$ induces a homomorphism $\Psi: G(\Delta) \to
\bng$.  It only remains to prove that $\Psi$ is also injective.  We do so
by showing that the composition map $\Phi \circ \Psi$ from right-angled
Artin groups to right-angled Artin groups is injective.

Assume that $\Phi \circ \Psi$ is not injective.  Fix a nontrivial element 
$\alpha$ in the kernel of $\Phi \circ \Psi$ and let $w$ be a geodesic word 
representing $\alpha$.  Since $\alpha$ is nontrivial, $w$ is non-empty, 
and $w$ contains a letter $a$ (or its inverse) representing a generator in 
$A$.  Then $w$ has the form 
  $$w = w_0a^{\epsilon_1}w_1a^{\epsilon_2}w_2\dots a^{\epsilon_k}w_k$$
where $\epsilon_i \in \{\pm 1\}$ for $i = 1, \dots, k$ and each $w_i$ is a
(possibly empty) word in the generators of $A$ not containing $a$ or
$a^{-1}$.

Let $\gamma$ be the Artin loop in $\Gamma$ corresponding to $a$.  Let $p$
be the point on $\gamma$ which is in the Artin basepoint $\{x_1, \dots,
x_n\}$.  Let $e_1e_2\dots e_m$ be the edges of $\gamma$ traversed in order
from $p$ - i.e. $e_1e_2\dots e_m$ is the image of $\gamma$ under $\Phi$.

Let $g_i \in A_{\Gamma}$ be the image under $\Phi \circ \Psi$ of the
element in $A$ represented by $w_i$.

Consider $A_{\Gamma}$ as an HNN-extension with stable letter $e_1$. 
Britton's lemma tells us that since $(\Phi \circ \Psi)(w) = 1$, there 
must exist an $e_1$-pinch in the appropriate expansion of
$$g_0(e_1e_2\dots e_m)^{2\epsilon_1}g_1 \dots (e_1e_2\dots 
  e_m)^{2\epsilon_k} g_k.$$
For any $e_1$-pinch of the form $e_1^{-1}g_ie_1$, replace the pinch with 
$g_i$.  Note each $g_i$ may be represented without using any of $e_1, 
\dots, e_m$.  Since $\Psi$ involves a square, not every occurrence of 
$e_1$ may be replaced in this way, so we know that there must 
still exist $e_1$-pinches in the result.  The only possible remaining 
$e_1$-pinches must have one of the following forms:
\begin{enumerate}
\item $e_1e_2\dots e_m g_i e_m^{-1}\dots e_2^{-1}e_1^{-1}$
\item $e_1^{-1}e_m^{-1}\dots e_2^{-1} g_i e_2e_3\dots e_me_1.$
\end{enumerate}
We consider the first case; the second case is dealt with in a similar 
way.

If $e_1e_2\dots e_m g_i e_m^{-1}\dots e_2^{-1}e_1^{-1}$ is an $e_1$-pinch,
then $e_2\dots e_m g_i e_m^{-1}\dots e_2^{-1} = g_i'$, where $g_i' \in
G(link(e_1))$.  Thus, $e_2\dots e_m g_i e_m^{-1}\dots e_2^{-1}(g_i')^{-1}
= 1$.  Since $e_2, e_m \not\in link(e_1)$, it follows that $e_2\dots e_m
g_i e_m^{-1}\dots e_2^{-1}$ is an $e_2$-pinch in $A_{\Gamma - e_1}$, by
Britton's Lemma.  Iterating, we see that $e_j\dots e_m g_i e_m^{-1}\dots
e_j^{-1}$ is an $e_j$-pinch for each $j = 1, \dots, m$.

In particular, we have proven that 
  $$e_1\dots e_mg_ie_m^{-1}\dots e_1^{-1} = g_i.$$
Thus, $(e_1\dots e_m)^2g_i(e_m^{-1}\dots e_1^{-1})^2 = g_i$.  
But then
$$a^{\epsilon_1}w_1a^{\epsilon_2}w_2\dots 
  a^{\epsilon_{i-1}}w_{i-1}w_iw_{i+1}a^{\epsilon_{i+2}}w_{i+2} \dots 
  a^{\epsilon_k}w_k$$
is a word representing $\alpha$ which is shorter than $w$.  As $w$ was
assumed to be a geodesic, this is a contradiction.  Thus no such $\alpha$ 
exists, and $\Phi \circ \Psi$ is injective.
\end{proof}

\begin{figure}[htbp]
\begin{center}
\input{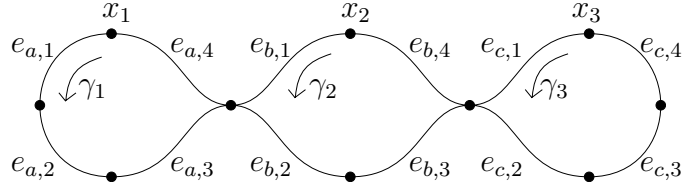}
\end{center}

\caption{An example showing why $\Psi$ maps to squares.}{Shown is a
possible graph $\Gamma$ for the right-angled Artin group $A = \langle
a,b,c | [a,c]\rangle$.  The braid group on 3 strands on $\Gamma$ is such
that if $\Psi$ is defined by sending each generator of $A$ to the
corresponding Artin loop \emph{without} squaring, then $\Psi$ is not
injective.}

\label{fig:noninjectiveexample}
\end{figure}

Note that if $\Psi$ is defined only by sending $a_i$ to $\gamma_i$ instead
of $\gamma_i^2$, injectivity may not hold.  For instance, consider Figure
\ref{fig:noninjectiveexample}.  This is a possible graph $\Gamma$ for the
right-angled Artin group $\langle a,b,c | [a,c]\rangle$.  If $a \mapsto
\gamma_a \mapsto e_{a,1}e_{a,2}e_{a,3}e_{a,4}$, $b\mapsto \gamma_b \mapsto
e_{b,1}e_{b,2}e_{b,3}e_{b,4}$, and $c \mapsto \gamma_c \mapsto
e_{c,1}e_{c,2}e_{c,3}e_{c,4}$, we leave it to the reader to verify that $g
= cbab^{-1}c^{-1}ba^{-1}b^{-1}$ is not trivial but $(\Phi \circ \Psi)(g)$
is (it is insightful to visualize exactly what is happening here in terms
of braids on $\Gamma$).

\begin{corollary}\label{cor:chromaticnumber}
Let $n_0$ be the chromatic number for $\Delta$.  Then there exists an 
embedding of $G(\Delta)$ into the graph braid group $B_{n_0}\Gamma$ for 
some $\Delta$-halo $\Gamma$.
\end{corollary}

\begin{proof}
Let $C_0$ be a coloring of $\Delta$ using exactly $n_0$ colors.  Let 
$\Gamma$ be a $\Delta$-halo for $C_0$.  The result then follows from 
Theorem \ref{prop:psiinjective}.
\end{proof}

\begin{example}\label{ex:C_6} As an application of Corollary
\ref{cor:chromaticnumber}, consider the right-angled Artin group $A$ for 
which $\Delta = C_6$, the cycle of length $6$. Thus,
  $$A = \langle a_1,a_2,a_3,a_4,a_5,a_6 | [a_i,a_{i+1}] (i = 1,\dots,5), 
  [a_6,a_1]\rangle.$$

\begin{figure}[htbp]
\begin{center}
\input{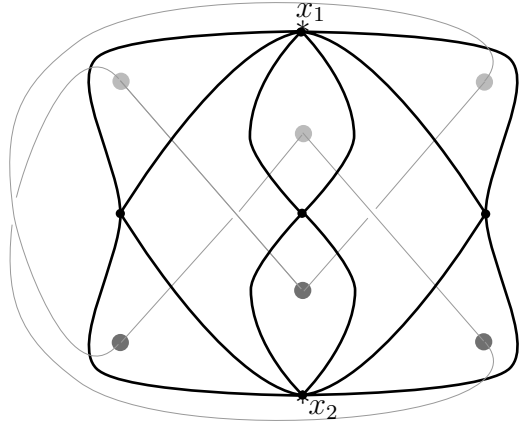}
\end{center}

\caption{A $\Delta$-halo for $C_6$.}
{This figure shows $(C_6)$ in grey, with a coloring given by $a_i
\mapsto 1$ for $i$ odd and $a_i \mapsto 2$ for $i$ even. A possible graph
$\Gamma = \Gamma(C_6)$ is shown.}
\label{fig:C_6}

\end{figure}

Figure \ref{fig:C_6} shows $C_6$ with a coloring using $2$ colors and a
possible graph $\Gamma(C_6)$.  Theorem \ref{thm:embedding} tells us that
we may embed $G(C_6)$ into $B_n\Gamma(C_6)$ where $n=2$!  Note that
$\Gamma(C_6)$ is planar.

Recall that a hyperbolic surface group is the fundamental group of a
compact hyperbolic surface without boundary.  Servatius et al.  
\cite{ServatiusDromsServatius} (see also \cite{GordonLongReid}) have shown
that a right-angled Artin group whose graph has an achordal cycle of
length at least 5 contains an embedded hyperbolic surface subgroup.  
Since $A$ embeds into $\btg(C_6)$, it follows that there is a hyperbolic
surface subgroup embedded in $\btg(C_6)$. 
\end{example}

Example \ref{ex:C_6} proves Theorem \ref{thm:planartwostrandsurface} via
the following corollary:

\begin{corollary}\label{cor:planartwostrandsurface}
There exists a \emph{planar} graph braid group - namely $\btg(C_6)$
- which contains a hyperbolic surface subgroup, with only $n = 2$ strands.
\end{corollary}

It is known that hyperbolic surface groups can be subgroups of graph braid
groups.  In fact, in \cite{Abrams1}, Abrams proved that

\begin{proposition}
The pure graph braid groups $PB_2K_{3,3}$ and $PB_2K_5$ are the
fundamental groups of closed orientable hyperbolic manifolds of genus $4$
and $6$, respectively.
\end{proposition}

This proposition, combined with embedding results in \cite{Abrams1}, yields
the existence of hyperbolic surface subgroups in almost every nonplanar
graph braid group (the possible exceptions being $B_nK_5$ and $B_nK_{3,3}$
for $n \geq 3$).  No planar examples were previously known.

Of course, Corollary \ref{cor:planartwostrandsurface} relied on the
observation that the graph $\Gamma(C_6)$ in the above example was planar.
It is not a priori true in general that a planar graph $\Delta$ gives rise 
to a planar $\Delta$-halo.  However, the following corollary does hold.  
Let $\Delta^{op}$ be the graph with vertex set equal to $\Delta$ and an 
edge connecting two vertices if and only if $\Delta$ does not have that 
edge.

\begin{corollary}
Let $A = G(\Delta)$ be a right-angled Artin group with associated graph
$\Delta$ such that $\Delta^{op}$ has a finite-sheeted planar covering.  
Then $A$ embeds in some \emph{planar} graph braid group $\bng$ for $n \leq 
|V(\Delta)|$.
\end{corollary}

\begin{proof} 
In Proposition 19 of \cite{CrispWiest}, Crisp and Wiest showed that if
$\Delta^{op}$ has a finite-sheeted planar covering $\tilde{\Delta}^{op}$,
then if $\tilde{\Delta}$ is the opposite graph of $\tilde{\Delta}^{op}$,
there exists an injective homomorphism $j: G(\Delta) \to
G(\tilde{\Delta})$.  As $\tilde{\Delta}$ is a cover of $\Delta$, 
$\tilde{\Delta}$ may be colored with $n \leq |V(\Delta)|$ colors.  We 
leave it to the reader to verify that a planar opposite graph gives rise 
to a planar halo.
\end{proof}

We end with a question.  As mentioned, the notion of a
$\Delta$-halo is similar in flavor to the proof of Theorem 10 of
\cite{CrispWiest}.  Although Crisp and Wiest deal with pure surface braid
groups instead of regular surface braid groups, is it possible to apply
the ideas presented here to reduce the number of strands needed for their
embedding?

\bibliography{refs-S2}
\bibliographystyle{plain}

\end{document}